\documentclass[12pt, reqno]{amsart}
\usepackage{amsmath, amsthm, amscd, amsfonts, amssymb, graphicx, color}
\usepackage[bookmarksnumbered, colorlinks, plainpages]{hyperref}
\hypersetup{colorlinks=true,linkcolor=red, anchorcolor=green, citecolor=cyan, urlcolor=red, filecolor=magenta, pdftoolbar=true}

\newtheorem{theorem}{Theorem}[section]

\theoremstyle{definition}
\newtheorem{definition}[theorem]{Definition}

\theoremstyle{remark}
\newtheorem{remark}[theorem]{Remark}
\numberwithin{equation}{section}

\begin{document}

\setcounter{page}{1}

\title[Well-posedness of Tricomi-Gellerstedt-Keldysh equations]{Well-posedness of Tricomi-Gellerstedt-Keldysh-type fractional elliptic problems}

\author[M. Ruzhansky, B. T. Torebek, B. Kh. Turmetov]{Michael Ruzhansky, Berikbol T. Torebek$^*$, Batirkhan Kh. Turmetov}

\address{\textcolor[rgb]{0.00,0.00,0.84}{Michael Ruzhansky \newline Department of Mathematics: Analysis,
Logic and Discrete Mathematics \newline Ghent University, Krijgslaan 281, Ghent, Belgium \newline
 and \newline School of Mathematical Sciences \newline Queen Mary University of London, United Kingdom}}
\email{\textcolor[rgb]{0.00,0.00,0.84}{michael.ruzhansky@ugent.be}}

\address{\textcolor[rgb]{0.00,0.00,0.84}{Berikbol T. Torebek \newline Department of Mathematics: Analysis,
Logic and Discrete Mathematics \newline Ghent University, Krijgslaan 281, Ghent, Belgium \newline and \newline Al--Farabi Kazakh National University \newline Al--Farabi ave. 71, 050040, Almaty, Kazakhstan \newline and \newline Institute of
Mathematics and Mathematical Modeling \newline 125 Pushkin str.,
050010 Almaty, Kazakhstan}}
\email{\textcolor[rgb]{0.00,0.00,0.84}{berikbol.torebek@ugent.be}}

\address{\textcolor[rgb]{0.00,0.00,0.84}{Batirkhan Kh. Turmetov \newline Department of Mathematics, Akhmet Yasawi University, \newline 29 B.Sattarkhanov str., 161200 Turkistan, Kazakhstan}}
\email{\textcolor[rgb]{0.00,0.00,0.84}{turmetovbh@mail.ru}}

\thanks{The first author was supported in parts by the FWO Odysseus Project 1 grant G.0H94.18N: Analysis and Partial Differential Equations, by the EPSRC grant EP/R003025/2 and by the Methusalem programme of the Ghent University Special Research
Fund (BOF) (Grant number 01M01021). The second author was supported in parts by the FWO Odysseus Project 1 grant G.0H94.18N: Analysis and Partial Differential Equations and by a grant No.AP08052046 from the Ministry of Science and Education of the Republic of Kazakhstan.}

\let\thefootnote\relax\footnote{$^{*}$Corresponding author}



\begin{abstract} In this paper Tricomi-Gellerstedt-Keldysh-type fractional elliptic equations are studied. The results on the well-posedness of fractional elliptic boundary value problems are obtained for general positive operators with discrete spectrum and for Fourier multipliers with positive symbols. As examples, we discuss results in half-cylinder, star-shaped graph, half-space and other domains.
\end{abstract}
\maketitle
\tableofcontents
\section{Introduction}
\subsection{Statement of the problem and historical background}
The main purpose of this paper is to study the following fractional elliptic equation
\begin{equation}\label{1.1} \mathcal{D}^{2\alpha } u(x,y) - x^{2\beta}\mathcal{L}u(x,y) = 0,\,\left({x,y}
\right) \in \mathbb{R}_+\times\Omega,\end{equation} where $1/2<\alpha\leq 1,\,\beta>-\alpha,$ $\Omega\subset\mathbb{R}^N$ is a bounded domain with smooth boundary or $\Omega=\mathbb{R}^N$, and $\mathcal{D}_x^{2\alpha }$ means $\mathcal{D}_x^{2\alpha }= \partial_{0+,x}^\alpha \partial_{0+,x}^\alpha.$ Here $\partial_{0+,x}^\alpha$ is a Caputo fractional derivatives of order $\alpha:$
$$\partial_{0+,x}^\alpha u(x,y) = \frac{1}{{\Gamma \left(1- \alpha \right)}}\int\limits_0^x {\left(
{x - s} \right)^{-\alpha} \partial_s u\left( s, y\right)} ds,$$ and $\mathcal L$ satisfies one of the following properties
\begin{description}
  \item[(A)] a linear self-adjoint positive operator with a discrete spectrum $\{\lambda_k\geq 0:k\in \mathbb N\}$ on the Hilbert space $L^2(\Omega)$. According to $\lambda_k$, the operator $\mathcal L$ has the system of orthonormal eigenfunctions $\{e_k:k\in \mathbb N\}$ on $L^2(\Omega)$.\\
      As an example of $\mathcal L$, we can consider all self-adjoint positive operators that were given in \cite{RuzT1, RuzT2}. For example:
      \begin{itemize}
        \item Dirichlet-Laplacian, Neumann-Laplacian or fractional Dirichlet-Laplacian in a bounded domain;
        \item Sturm-Liouville operator or its involution perturbations in a finite interval;
        \item integro-differential operators with fractional derivatives.
      \end{itemize}
  \item[(B)] Fourier multiplier $a(D)$ with symbol $a(\xi)\geq 0,\,\,\xi\in \mathbb{R}^N,$ i.e. $a(D)=\mathcal{F}^{-1}\left(a(\xi)\mathcal{F}\right),\,\,\xi\in \mathbb{R}^N,$ where $\mathcal{F}$ is the Fourier transform and $\mathcal{F}^{-1}$ is the inverse Fourier transform.\\
      As an example of $\mathcal L$, we can consider all operators with nonnegative symbol (see \cite{Ruzh}). For example:
      \begin{itemize}
        \item Laplace operator $-\Delta$ with symbol $|\xi|^2$ or fractional Laplacian $(-\Delta)^s,\,s\in(0,1),$ with symbol $|\xi|^{2s}$;
        \item Linear partial differential operator $\sum\limits_{|\beta|\leq m} a_\beta D^\beta,\,\,a_\beta\geq 0,$ with nonnegative symbol $\sum\limits_{|\beta|\leq m} a_\beta \xi^\beta\geq 0,$ with $D^\beta=\left(\frac{1}{i}\partial_{x_1}\right)^{\beta_1}\cdot ... \cdot \left(\frac{1}{i}\partial_{x_N}\right)^{\beta_N}$.
      \end{itemize}
\end{description}

The need to study the boundary value problems for the fractional elliptic equations to describe the production processes in
mathematical modeling of socio-economic systems was shown in \cite{Nak}. In \cite{Nak} the attention was drawn to the fact that the problem of finding a generalized two-factor Cobb-Douglas function is reduced to the Dirichlet problem for the fractional elliptic equation.

The equation \eqref{1.1} is a generalization of the following well-known equations:
\begin{itemize}
  \item If $\alpha=1,$ $\beta=0$ and $\mathcal{L}=-\Delta=-\sum\limits_{j=1}^n\frac{\partial^2}{\partial y^2_j},$ then the equation \eqref{1.1} coincides with the
classical Laplace equation $$u_{xx}(x,y)+\Delta_yu(x,y)=0,\,x>0,\,y\in\mathbb{R}^N;$$
  \item If $N=1,$ $\alpha=1,$ $\beta=\frac{1}{2}$ and $\mathcal{L}=-\frac{\partial^2}{\partial y^2},$ then the equation \eqref{1.1} coincides with the
classical Tricomi equation (\cite{Tricomi}) $$u_{xx}(x,y)+xu_{yy}(x,y)=0,\,x>0,\,y\in\mathbb{R};$$
  \item If $N=1,$ $\alpha=1,$ $\beta=m>0$ and $\mathcal{L}=-\frac{\partial^2}{\partial y^2},$ then the equation \eqref{1.1} coincides with the
classical Gellerstedt equation (\cite{Geller}) $$u_{xx}(x,y)+x^{m}u_{yy}(x,y)=0,\,x>0,\,y\in\mathbb{R};$$
  \item If $N=1,$ $\alpha=1,$ $\beta=-k\in(-2,0)$ and $\mathcal{L}=-\frac{\partial^2}{\partial y^2},$ then the equation \eqref{1.1} coincides with the
classical Keldysh equation (\cite{Keldysh}) $$u_{xx}(x,y)+x^{-k}u_{yy}(x,y)=0,\,x>0,\,y\in\mathbb{R}.$$
\end{itemize}
The above equations are used in transonic gas dynamics \cite{Bers58}, and in mathematical models of cold plasma \cite{Otway}.

Note that the study of Tricomi, Gellerstedt and Keldysh equations was done in many papers \cite{Alg, Gelf1, Gelf2, Gelf3, Mois, Xu}. The boundary value problems for the fractional elliptic equations are studied in \cite{Amb, Caff, Mas, Turm}.

\subsection{Three-parameter Mittag-Leffler (Kilbas-Saigo) function}

First, we recall the definition of the Kilbas-Saigo function (three-parameter Mittag-Leffler function) and some of its particular cases.
\begin{itemize}
\item {\bf Classical Mittag-Leffler function.} The classical Mittag-Leffler function $E_{\alpha,1}(z)$ defined by (\cite{M-L03})
\begin{equation*}
E_{\alpha,1}(z)=\sum\limits_{k=0}^\infty \frac{z^k}{\Gamma(\alpha k+1)},\,\,\alpha>0,\, z\in \mathbb{C},
\end{equation*} is a natural extension of the exponential function $E_{1,1}(z)=\exp(z),$ and also of the hyperbolic cosine function $E_{2,1}(z)=\cosh{\sqrt{z}}.$

The most interesting properties of Mittag-Leffler function are associated with its upper-lower estimates for $0<\alpha<1$ as follows (\cite{TSim15}):
\begin{equation}\label{MLF1}
\frac{1}{1+\Gamma(1-\alpha)z}\leq E_{\alpha, 1}(-z)\leq \frac{1}{1+\frac{1}{\Gamma(1+\alpha)}z},\, z\geq 0.
\end{equation}
\end{itemize}
\begin{itemize}
\item {\bf Two-parameter Mittag-Leffler function.} The two-parameter Mittag-Leffler function $E_{\alpha,\beta}(z)$ is defined by
\begin{equation*}
E_{\alpha,\beta}(z)=\sum\limits_{k=0}^\infty \frac{z^k}{\Gamma(\alpha k+\beta)},\,\,\alpha>0,\, \beta>0,\, z\in \mathbb{C}.
\end{equation*}
This function, sometimes called a Mittag-Leffler-type function, first appeared in \cite{W05}. When $\beta=1$, $E_{\alpha,\beta}(z)$ coincides with the classical Mittag-Leffler function $E_{\alpha,1}(z).$
\end{itemize}
\begin{itemize}
  \item {\bf Three-parameter (Kilbas-Saigo) Mittag-Leffler function.}
Another generalization of the Mittag-Leffler function was introduced by Kilbas and Saigo \cite{KS95} in terms of a special function of the form
\begin{equation}
\label{SF-01}
E_{\alpha, m, n}(z)=1+\sum_{k=1}^{\infty}\prod_{j=0}^{k-1}\frac{\Gamma(\alpha(jm+n)+1)}{\Gamma(\alpha(jm+n+1)+1)}\,z^{k},
\end{equation}
where $\alpha,\, m$ are real numbers and $n\in \mathbb{C}$ such that
\begin{equation}\label{cond1}
\alpha>0,\, m>0,\, \alpha(jm+n)+1\neq -1,-2,-3,... (j\in\mathbb N_0).
\end{equation}
In particular, if $m = 1,$ the function $E_{\alpha, m, n}(z)$ is reduced to the two-parameter Mittag-Leffler function:
\begin{equation*}
 E_{\alpha, 1, n}(z)=\Gamma(\alpha n+1)E_{\alpha, \alpha n+1}(z),
\end{equation*} and if $m = 1, n=0,$ then it coincides with the classical Mittag-Leffler function:
\begin{equation*}
 E_{\alpha, 1, 0}(z)=E_{\alpha, 1}(z).
\end{equation*}
Recently Simon et al. \cite{TSim19} obtained the following interesting estimates of the Kilbas-Saigo functions:
\begin{equation}\label{estim1}
\frac{1}{1+\Gamma(1-\alpha)z}\leq E_{\alpha, m, m-1}(-z)\leq \frac{1}{1+\frac{\Gamma(1+(m-1)\alpha)}{\Gamma(1+m\alpha)}z},\, z\geq 0,
\end{equation} where $m>0$ and $0<\alpha<1$.
\end{itemize}

\subsection{Ill-posedness of the non-sequential problem}
As generally $$\partial_x^\alpha \partial_x^\alpha\neq \partial_x^{2\alpha},$$ the equation \eqref{1.1} is different from the following non-sequential equation
\begin{equation}\label{001} \partial_x^{2\alpha}u(x,y) - x^{2\beta}\mathcal{L}u(x,y) = 0,\,\left({x,y}
\right) \in \mathbb{R}_+\times\Omega.\end{equation} However, we cannot consider the problem of bounded solutions of equation \eqref{001} in $x\in \mathbb{R}_+$, since for such class of functions, nontrivial solutions of equation \eqref{001} may not exist.
We demonstrate this with the following example:\\
Let $1<2\alpha<2,$ $\beta=0,$ and $\mathcal{L}=-\Delta=-\sum\limits_{j=1}^n\frac{\partial^2}{\partial y^2_j}$ in \eqref{001}. Then using the Fourier transform to \eqref{001}  with respect to $y$ we have
\begin{equation}\label{002} \partial_x^{2\alpha}\hat{u}(x,\xi) - |\xi|^{2}\hat{u}(x,\xi) = 0,\,x>0,\,\xi\in\mathbb{R}^N.\end{equation}
The general solution to the equation \eqref{002} has the form \cite[Example 4.10]{1}
\begin{equation*}\hat{u}(x,\xi) = C_1(\xi) E_{\alpha, 1} \left( |\xi|^{2} x^{\alpha}\right) + C_2(\xi) t E_{\alpha, 2} \left( {|\xi|^{2} x^{\alpha}}\right),\end{equation*} where $C_1(\xi),$ $C_2(\xi)$ are arbitrary constants and $E_{\alpha,\beta}(z)$ is the Mittag-Leffler function. From the asymptotic estimate of the Mittag-Leffler function $$E_{\alpha,\beta}(z)\sim z^{\frac{1-\beta}{\alpha}}e^{z^{\frac{1}{\alpha}}},\,z\rightarrow\infty,$$ it follows that
$$\lim\limits_{x\rightarrow\infty}E_{\alpha, 1} \left( |\xi|^{2s} x^{\alpha}\right)\rightarrow\infty\,\,\,\textrm{and}\,\,\, \lim\limits_{x\rightarrow\infty}E_{\alpha, 2} \left( |\xi|^{2s} x^{\alpha}\right)\rightarrow\infty.$$
Therefore, the equation \eqref{002} does not have a bounded solution in $x\in\mathbb{R}_+.$

\subsection{One dimensional fractional differential equation}
Let $ 0 < \alpha  \le 1 ,$ $ \mu $ is a positive real number. For
further exposition we need to give some information about the exact solutions of differential equations of the form: \begin{equation}\label{2.1} \mathcal{D}^{2\alpha }
h\left( x \right) - \mu^2 x^{2\beta} h\left( x \right) =
0,\,x > 0.\end{equation}

Using the method of constructing the solution of the fractional-order differential equations developed in \cite{Tur1, Tur2}, one can show that the functions \begin{equation}\label{2.3} \left\{E_{\alpha, 1+\frac{\beta}{\alpha}, \frac{\beta}{\alpha}} \left( { \mu x^{\alpha+\beta}  }\right), E_{\alpha, 1+\frac{\beta}{\alpha}, \frac{\beta}{\alpha}} \left( { -\mu x^{\alpha+\beta}  }
\right)\right\},\end{equation} are solutions of the equation \eqref{2.1}.

It is easy to show that the functions \eqref{2.3} are linearly independent. Hence, the system of functions
\eqref{2.3} is a fundamental system for the equation \eqref{2.1}, and
therefore the general solution of this equation has the form:
\begin{equation}\label{2.4}h\left(x\right) = C_1 E_{\alpha, 1+\frac{\beta}{\alpha}, \frac{\beta}{\alpha}} \left( { \mu x^{\alpha+\beta}  }\right) + C_2 E_{\alpha, 1+\frac{\beta}{\alpha}, \frac{\beta}{\alpha}} \left( {- \mu x^{\alpha+\beta}  }\right),\end{equation} where $C_1$ and $C_2 $ are arbitrary constants.

It is easy to see that, if $ x \to +\infty ,$ then $$E_{\alpha, 1+\frac{\beta}{\alpha}, \frac{\beta}{\alpha}} \left( { \mu x^{\alpha+\beta}  }\right)\to +\infty,$$ since \begin{equation}\label{2.5}E_{\alpha, 1+\frac{\beta}{\alpha}, \frac{\beta}{\alpha}} \left( { \mu x^{\alpha+\beta}  }\right)\geq \frac{\mu\Gamma(\beta+1)}{\Gamma(\alpha+\beta+1)}x^{\alpha+\beta},\,x>0.\end{equation}
And for the function $E_{\alpha, 1+\frac{\beta}{\alpha}, \frac{\beta}{\alpha}} \left( {- \mu x^{\alpha+\beta}  }\right),$ the following estimate holds (\cite{TSim19}):
\begin{equation}\label{2.6}E_{\alpha, 1+\frac{\beta}{\alpha}, \frac{\beta}{\alpha}} \left( {- \mu x^{\alpha+\beta}  }\right)\leq \frac{1}{1+\frac{\Gamma(\beta+1)}{\Gamma(\alpha+\beta+1)}\mu x^{\alpha+\beta}},\, x>0.\end{equation}

\section{Well-posedness in a bounded domain}
Let $\mathcal{L}$ be a self-adjoint, positive operator with the discrete spectrum $\{\lambda_k\geq 0:\;k\in\mathbb N\}$ on $L^2(\Omega)$. The main assumption in this section is that the system of eigenfunctions $\{e_k\in L^2(\Omega):k\in\mathbb N\}$ of the operator $\mathcal{L}$ is an orthonormal basis in $L^2(\Omega)$.

The Hilbert space $\mathcal{H}^\mathcal{L}(\Omega)$ is defined by
$$\mathcal{H}^\mathcal{L}(\Omega)=\{u\in L^2(\Omega):\, \sum\limits_{k = 0}^\infty \lambda^2_k|(u,e_k)|^2<\infty\},$$
with the norm
$$\|u\|^2_{\mathcal{H}^\mathcal{L}(\Omega)}=\sum\limits_{k = 0}^\infty \lambda^2_k|(u,e_k)|^2.$$

\begin{definition} The generalised solution of equation \eqref{1.1} in $\Omega\subset\mathbb{R}^N$ is a bounded function $u\in C\left(\mathbb{R}_+;L^2(\Omega)\right),$ such that $x^{-2\beta}\mathcal{D}_x^{2\alpha } u, \mathcal{L}u\in C\left(\mathbb{R}_+;L^2(\Omega)\right).$\end{definition}
\begin{theorem}\label{th1}Let $\phi \in \mathcal{H}^\mathcal{L}(\Omega).$ Then the generalised solution of equation \eqref{1.1} satisfying conditions
\begin{equation}\label{1.2}
u(0,y)=\phi(y),\,y\in\Omega,
\end{equation}
and
\begin{equation}\label{1.3}
\lim\limits_{x\rightarrow+\infty}u(x,y)\,\,\,\,\text{is bounded for almost every}\,\,\,\,y\in\Omega,
\end{equation}
exists, it is unique and can be represented as \begin{equation}\label{1.5}u\left( {x,y} \right) = \sum\limits_{k = 0}^\infty {\phi_k E_{\alpha, 1+\frac{\beta}{\alpha}, \frac{\beta}{\alpha}} \left( {- \sqrt{\lambda_k} x^{\alpha+\beta}  }\right) e_k\left( y \right)},\,(x,y)\in [0,\infty)\times\Omega,\end{equation} where $\phi_k=\int\limits_\Omega {\phi\left( y \right) \overline{e_k \left( y \right)}dy},\,k\in\mathbb{Z}_+=0, 1, 2, ... ,$ and $E_{\alpha, m, l} \left( z \right)$ is a Kilbas-Saigo function.

In addition, the solution $u$ satisfies the following estimates:
\begin{align*}\|u\|_{C(\mathbb{R}_+;L^2(\Omega))}\leq\|\phi\|_{L^2(\Omega)},\end{align*}
\begin{align*}\sup\limits_{x\in(0,\infty)}\left\|x^{-2\beta}\mathcal{D}_x^{2\alpha}u\left( {x, \cdot} \right)\right\|_{L^2(\Omega)} \leq \|\phi\|_{\mathcal{H}^\mathcal{L}(\Omega)},\end{align*} and \begin{align*}\sup\limits_{x\in(0,\infty)}\|\mathcal{L}u\left({x, \cdot} \right)\|_{L^2(\Omega)}\leq\|\phi\|_{\mathcal{H}^\mathcal{L}(\Omega)}.\end{align*}
\end{theorem}

\begin{remark} If in Theorem \ref{th1} we replace the boundedness condition \eqref{1.3} by condition
\begin{equation}\label{1.4}
\lim\limits_{x\rightarrow+\infty}u(x,y)=0,\,\,\,\,y\in\Omega,
\end{equation}
then the problem \eqref{1.1}, \eqref{1.2}, \eqref{1.4} for the self-adjoint operators $\mathcal{L}$ with nonnegative eigenvalues $\lambda_k\geq 0,\, k\in \mathbb{N},$ becomes ill-posed. Indeed, it is easy to show that the bounded solution to Problem \eqref{1.1}, \eqref{1.2} has the form \eqref{1.5}. However, if we take into account condition \eqref{1.4}, then, for the existence of a solution to problem \eqref{1.1}, \eqref{1.2}, \eqref{1.4}, it is necessary and sufficient to have the condition
$$\int\limits_\Omega {\phi\left( y \right) dy}=0.$$
\end{remark}
\subsection{Particular cases} We now specify Theorem \ref{th1} to several concrete cases.
\subsubsection{Laplace equation in the half-strip and in the star-shaped graphs}
Our first example will focus on the Laplace equation.

$\bullet$ Let $\Omega=(0,1),$ $\alpha=1,$ $\beta=0$ and $$\mathcal{L}=-\frac{\partial^2}{\partial y^2},\,D(\mathcal{L}):=\{u\in W^1_2([0,1]),\, u(0)=u(1)=0\}.$$ Then the equation \eqref{1.1} coincides with the classical Laplace equation on the half-strip
\begin{equation}\label{1.1a} u_{xx}(x,y) + u_{yy}(x,y) = 0,\,\left({x,y}
\right) \in \mathbb{R}_+\times(0,1).\end{equation} It is known that the unique solution to problem \eqref{1.1a}, \eqref{1.2}, \eqref{1.3} is represented in the form
\begin{equation*}u\left( {x,y} \right) = \sum\limits_{k = 1}^\infty \phi_k e^{-k\pi x} \sin k\pi y.\end{equation*}

$\bullet$ Let $\Omega$ be a star-shaped metric graph consisting of $d$ segments of equal length, $\alpha=1,$ $\beta=0,$ and let $\mathcal{L}$ be a differential operator $\mathcal{L}=-\frac{\partial^2 v_j(y)}{\partial y^2},\,j=1,...,d,$ with boundary conditions \begin{align*}&v_j(0)=0,\,j=1,...,d,\\& v_1(\pi)=v_2(\pi)=...=v_d(\pi),\\& v'_1(\pi)+v'_2(\pi)+...+v'_d(\pi)=0.\end{align*} It is known (\cite{yang}) that the above operator is self-adjoint in $L^2_d([0,\pi])=\bigotimes\limits_{i=1}^dL^2([0,\pi])$ and has discrete spectrum $\lambda_k^d=\left(k-\frac{1}{2}\right)^2,\,k\in\mathbb{N}.$ Then the equation \eqref{1.1} coincides with the Laplace equation on the star-shaped graphs
\begin{equation}\label{1.1aa}
\Delta u(x,y)\equiv \Delta \left(\begin{array}{l}u_1(x,y)\\ u_2(x,y)\\ \vdots\\ u_d(x,y)\end{array}\right) = 0.\end{equation} Then the unique solution to problem \eqref{1.1aa}, \eqref{1.2}, \eqref{1.3} is represented in the form
\begin{equation*}u(x,y)\equiv \left(\begin{array}{l}u_1(x,y)\\ u_2(x,y)\\ \vdots\\ u_d(x,y)\end{array}\right) = \sum\limits_{k = 1}^\infty \phi_k e^{- \left(k-\frac{1}{2}\right)x} \left(\begin{array}{l}1\\ 1\\ \vdots\\ 1\end{array}\right) \sin \left(k-\frac{1}{2}\right) y.\end{equation*}
\subsubsection{Fractional analogue of the Laplace equation with involution}
Let $\Omega=(-\pi,\pi),$ $\beta=0,$ and $$\mathcal{L}u(x)=-\frac{\partial^2}{\partial y^2}u(x)+\varepsilon \frac{\partial^2}{\partial y^2}u(-x),\,|\varepsilon|<1,$$ $$D(\mathcal{L}):=\{u\in W^1_2([-\pi,\pi]),\, u(-\pi)=u(\pi)=0\}.$$ Then the equation \eqref{1.1} coincides with the fractional analogue of the Laplace equation with involution on the half-strip
\begin{equation}\label{1.1b} \mathcal{D}_x^{2\alpha }u(x,y) + u_{yy}(x,y) - \varepsilon u_{yy}(x,-y) = 0,\,\left({x,y}
\right) \in \mathbb{R}_+\times(-\pi,\pi).\end{equation} It is known (\cite{Turm}) that there exist a unique solution to problem \eqref{1.1b}, \eqref{1.2}, \eqref{1.3} and it can be represented in the form
\begin{equation*}u\left( {x,y} \right) = \sum\limits_{k = 1}^\infty \phi_k E_{\alpha,1}\left({-\left(1+(-1)^k\varepsilon\right)k\pi x^\alpha}\right) \sin k\pi y.\end{equation*}
\subsubsection{Elliptic Tricomi and Gellerstedt equation}
Let $\alpha=1,$ $\beta>-2,$ and $\mathcal{L}=-\frac{\partial^2}{\partial y^2},\,D(\mathcal{L}):=\{u\in W^1_2([0,1]),\, u(0)=u(1)=0\}.$

$\bullet$ If $\beta=1$ then the equation \eqref{1.1} coincides with the classical Tricomi equation \begin{equation}\label{1.1c}u_{xx}(x,y)+x u_{yy}(x,y)=0,\,x>0,\,y\in(0,1),\end{equation}
and the unique solution to problem \eqref{1.1c}, \eqref{1.2}, \eqref{1.3} can be written as
\begin{equation*}u\left( {x,y} \right) = \sum\limits_{k = 1}^\infty \phi_k\, \text{Ai}({-k\pi x}) \sin k\pi y,\end{equation*} where $\text{Ai}(z)$ is the Airy function.

$\bullet$ If $\beta>-2$ then the equation \eqref{1.1} coincides with the classical Gellerstedt equation \begin{equation}\label{1.1cc}u_{xx}(x,y)+x^\beta u_{yy}(x,y)=0,\,x>0,\,y\in(0,1),\end{equation}
and the unique solution to problem \eqref{1.1cc}, \eqref{1.2}, \eqref{1.3} can be written as (see \cite{Mois})
\begin{equation*}u\left( {x,y} \right) = \sum\limits_{k = 1}^\infty \phi_k\, \sqrt{x}K_{\frac{1}{\beta+2}}\left(\frac{2\pi k x^{\frac{2}{\beta+2}}}{\beta+2}\right) \sin k\pi y,\end{equation*} where $K_{\nu}(z)$ is the Macdonald function.

\subsubsection{Fractional elliptic equation with variable coefficients}
If $\beta=0$ and $$\mathcal{L}=(1-y)^\mu (1+y)^\mu D^\mu_{1-,y}\partial^\mu_{-1+,y},$$ $$u(-1)=I^{1-\mu}_{1-,y}\partial^\mu_{-1+,y}u(1)=0,$$ then the equation \eqref{1.1} coincides with the equation \begin{equation}\label{1.1d}u_{xx}(x,y)+(1-y)^\mu (1+y)^\mu D^\mu_{1-,y}\partial^\mu_{-1+,y}u(x,y)=0,\,x>0,\,y\in(-1,1),\end{equation} where $\mu\in(0,1),$ $D_{1-,y}^\mu$ is a right-side Riemann-Liouville fractional derivative of order $\mu\in(0,1)$
$$D_{1-,y}^\mu u(x,y) = \frac{1}{{\Gamma \left(1- \mu \right)}}\frac{\partial}{\partial y}\int\limits_{y}^1 {\left(
{s - y} \right)^{-\mu} u\left( x, s\right)} ds,$$ $\partial_{-1+,y}^\mu$ is a left-side Caputo fractional derivative of order $\mu\in(0,1)$
$$\partial_{-1+,y}^\mu u(x,y) = \frac{1}{{\Gamma \left(1- \mu \right)}}\int\limits_{-1}^y {\left(
{y - s} \right)^{-\mu} u_s\left( x, s\right)} ds,$$ $I_{1-,y}^{1-\mu}$ is a right-side Riemann-Liouville fractional integral of order $\mu\in(0,1)$
$$I_{1-,y}^{1-\mu} u(x,y) = \frac{1}{{\Gamma \left(1- \mu \right)}}\int\limits_{y}^1 {\left(
{s - y} \right)^{-\mu} u\left( x, s\right)} ds.$$
The unique solution of problem \eqref{1.1d}, \eqref{1.2}, \eqref{1.3} can be written as
\begin{equation*}u\left( {x,y} \right) = \sum\limits_{k = 1}^\infty \phi_k\, \exp\left({-\frac{\Gamma(k+\mu)}{\Gamma(k-\mu)} x}\right) (1+y)^\mu P^{-\mu,\mu}_{k-1}(y),\end{equation*} where $P^{-\mu,\mu}_{k-1}(y)$ is the Jacobi polynomial (\cite{JCP}) $$P^{-\mu,\mu}_{k-1}(y)=\sum\limits_{n=0}^{k-1}\left(\begin{array}{l}k-1-\mu\\k-1-n\end{array}\right) \left(\begin{array}{l}k-1+\mu\\n\end{array}\right)\left(\frac{y-1}{2}\right)^n\left(\frac{y+1}{2}\right)^{k-1-n}.$$

\subsection{Proof of Theorem \ref{th1}}
\subsubsection{Existence of solution.}
As $\mathcal{L}$ is self-adjoint in $L^2(\Omega),$ any solution of problem \eqref{1.1}, \eqref{1.2}--\eqref{1.3} can be represented as: \begin{equation}\label{3.2}u\left( {x,y} \right) = \sum\limits_{k = 0}^\infty {u_k
\left( x \right) e_k\left( y \right)},\,(x,y)\in \mathbb{R}_+\times\Omega.\end{equation} It is clear that if $\phi\in \mathcal{H}^\mathcal{L}(\Omega)$, then
it can be represented in the form \begin{equation}\phi\left( y \right) = \sum\limits_{k = 0}^\infty {\phi_k
e_k \left( y \right)},\,y\in\Omega,\end{equation} where
$\phi_k=\int\limits_\Omega {\phi\left( y \right) \overline{e_k \left( y \right)}dy}.$

Substituting function \eqref{3.2} into equation \eqref{1.1}, we obtain the following problem for
$u_k(x),$
\begin{equation}\label{3.3}\mathcal{D}^{2\alpha } u_k
\left( x \right) - \lambda_k x^{2\beta}u_k \left( x \right) = 0,\, x >0,\end{equation} \begin{equation}\label{3.4}u_k \left(0\right) = \phi_k,\, u_k \left( \infty \right)\leq C,\,C=const,\end{equation} where $\lambda_k>0$ are eigenvalues of $\mathcal{L}$.

According to formula \eqref{2.4}, the general solution to equation \eqref{3.3} has the form: \begin{equation*}u_k \left(
x \right) = C_1 E_{\alpha, 1+\frac{\beta}{\alpha}, \frac{\beta}{\alpha}} \left( { \sqrt{\lambda_k} x^{\alpha+\beta}  }\right) + C_2 E_{\alpha, 1+\frac{\beta}{\alpha}, \frac{\beta}{\alpha}} \left( {- \sqrt{\lambda_k} x^{\alpha+\beta}  }\right),\end{equation*} where $C_1$ and $C_2 $ are arbitrary constants.

Since $$E_{\alpha, 1+\frac{\beta}{\alpha}, \frac{\beta}{\alpha}} \left( { \sqrt{\lambda_k} x^{\alpha+\beta}  }\right)\to +\infty,\,\,\, \text{as} \,\,\,x \to +\infty,$$ we have $C_1=0.$

Since $$E_{\alpha, 1+\frac{\beta}{\alpha}, \frac{\beta}{\alpha}} \left( { -\sqrt{\lambda_k} x^{\alpha+\beta}  }\right)\to 0,\,\,\, \text{as} \,\,\,x \to +\infty,$$ then by \eqref{3.4} we have
\begin{equation}\label{3.6}u_k \left(
x \right) = \phi_k E_{\alpha, 1+\frac{\beta}{\alpha}, \frac{\beta}{\alpha}} \left( {- \sqrt{\lambda_k} x^{\alpha+\beta}  }\right),\end{equation}
hence
\begin{equation*}u\left( {x,y} \right) = \sum\limits_{k = 0}^\infty {\phi_k E_{\alpha, 1+\frac{\beta}{\alpha}, \frac{\beta}{\alpha}} \left( {- \sqrt{\lambda_k} x^{\alpha+\beta}  }\right) e_k\left( y \right)},\,(x,y)\in \mathbb{R}_+\times\Omega.\end{equation*}

\subsubsection{Convergence of solution.}
The estimate \eqref{2.6} gives
\begin{align*}|u_k \left(
x \right)| \leq \frac{|\phi_k|}{1+\frac{\Gamma(\beta+1)}{\Gamma(\alpha+\beta+1)}\sqrt{\lambda_k} x^{\alpha+\beta}},\end{align*}
which implies
\begin{align*}\sup\limits_{x\geq 0}\|u\left( {x, \cdot} \right)\|_{L^2(\Omega)}^2 &\leq \sup\limits_{x\geq 0}\sum\limits_{k = 0}^\infty {|\phi_k|^2 \left|E_{\alpha, 1+\frac{\beta}{\alpha}, \frac{\beta}{\alpha}} \left( {- \sqrt{\lambda_k} x^{\alpha+\beta}  }\right)\right|^2 \|e_k\|_{L^2(\Omega)}^2}\\& \leq \sup\limits_{x\geq 0}\sum\limits_{k = 0}^\infty \frac{|\phi_k|^2}{\left(1+\frac{\Gamma(\beta+1)}{\Gamma(\alpha+\beta+1)}\sqrt{\lambda_k} x^{\alpha+\beta}\right)^2} \\&\leq \sum\limits_{k = 0}^\infty {|\phi_k|^2}=\|\phi\|_{L^2(\Omega)}^2<\infty,\end{align*} thanks to Parseval's identity.
Let us calculate $\mathcal{D}_x^{2\alpha } u$ and $\mathcal{L}u.$ We have
\begin{align*}\mathcal{D}_x^{2\alpha}u\left( {x,y} \right)& = \sum\limits_{k = 0}^\infty {\phi_k \mathcal{D}_x^{2\alpha}E_{\alpha, 1+\frac{\beta}{\alpha}, \frac{\beta}{\alpha}} \left( {- \sqrt{\lambda_k} x^{\alpha+\beta}  }\right) e_k\left( y \right)}\\&= x^{2\beta}\sum\limits_{k = 0}^\infty \lambda_k{\phi_k E_{\alpha, 1+\frac{\beta}{\alpha}, \frac{\beta}{\alpha}} \left( {- \sqrt{\lambda_k} x^{\alpha+\beta}  }\right) e_k\left( y \right)},\,(x,y)\in \mathbb{R}_+\times\Omega,\end{align*} and \begin{align*}\mathcal{L}u\left( {x,y} \right)& = \sum\limits_{k = 0}^\infty {\phi_k E_{\alpha, 1+\frac{\beta}{\alpha}, \frac{\beta}{\alpha}} \left( {- \sqrt{\lambda_k} x^{\alpha+\beta}  }\right) \mathcal{L}e_k\left( y \right)}\\&= \sum\limits_{k = 0}^\infty \lambda_k{\phi_k E_{\alpha, 1+\frac{\beta}{\alpha}, \frac{\beta}{\alpha}} \left( {- \sqrt{\lambda_k} x^{\alpha+\beta}  }\right) e_k\left( y \right)},\,(x,y)\in \mathbb{R}_+\times\Omega.\end{align*}
Applying the above calculations and Parseval's identity we have
\begin{align*}\sup\limits_{x\in(0,\infty)}\left\|x^{-2\beta}\mathcal{D}_x^{2\alpha}u\left( {x, \cdot} \right)\right\|^2_{L^2(\Omega)} \leq \sum\limits_{k = 0}^\infty \lambda_k^2{|\phi_k|^2 }=\|\phi\|^2_{\mathcal{H}^\mathcal{L}(\Omega)}<\infty,\end{align*} and \begin{align*}\sup\limits_{x\in(0,\infty)}\|\mathcal{L}u\left({x, \cdot} \right)\|^2_{L^2(\Omega)}\leq\sum\limits_{k = 0}^\infty \lambda_k^2{|\phi_k|^2}=\|\phi\|^2_{\mathcal{H}^\mathcal{L}(\Omega)}<\infty.\end{align*}

\subsubsection{Uniqueness of solution.} Suppose that there are two solutions $u_1(x, y)$ and $u_2(x,y)$ of problem \eqref{1.1}, \eqref{1.2}--\eqref{1.3}. Let $$u(x,y)=u_1(x,y)-u_2(x,y).$$ Then $u(x,y)$ satisfies the equation \eqref{1.1} and homogeneous conditions \eqref{1.2}--\eqref{1.3}.

Let us consider the function \begin{equation}\label{3.5}u_k(x)=\int\limits_\Omega u(x,y)\overline{e_k(y)}dy,\,k\in\mathbb{Z}_+,\,x\geq 0.\end{equation}
Applying $\mathcal{D}^{2\alpha}$ to the function \eqref{3.5} by \eqref{1.1} we have
\begin{align*}\mathcal{D}^{2\alpha}u_k(x)&=\int\limits_\Omega \mathcal{D}^{2\alpha}_xu(x,y)\overline{e_k(y)}dy =x^{2\beta}\int\limits_\Omega \mathcal{L}u(x,y)\overline{e_k(y)}dy \\& =x^{2\beta}\int\limits_\Omega u(x,y)\mathcal{L}\overline{e_k(y)}dy=x^{2\beta}\lambda_k\int\limits_\Omega u(x,y)\overline{e_k(y)}dy\\&= x^{2\beta}\lambda_k u_k(x),\,k\in\mathbb{Z}_+,\,x\geq 0.\end{align*} Also from \eqref{1.2} and \eqref{1.3} we have
$u_k(0)=0,\,\,u_k(\infty)\,\,\,\text{is bounded}.$ Then from \eqref{3.6} we conclude that $u_k(x)=0,\,x\geq 0.$ This implies $\int\limits_\Omega u(x,y)\overline{e_k(y)}dy=0,$ and the completeness of the system $e_k(x),\,k\in \mathbb{Z}_+,$ gives $u(x,y)\equiv 0,\,(x,y)\in [0,\infty)\times\Omega.$

\section{Well-posedness in $\mathbb{R}^N$} The Sobolev space $\mathcal{H}^\mathcal{L}(\mathbb{R}^N)$ is defined by
$$\mathcal{H}^\mathcal{L}(\mathbb{R}^N)=\{f\in L^2(\mathbb{R}^N):\, a(\xi)\hat{f}\in L^2(\mathbb{R}^N)\},$$ where $\hat{f}(\xi)=\frac{1}{(2\pi)^N}\int\limits_{\mathbb{R}^N}e^{-iy\xi} f(y)dy,\,\,\xi\in\mathbb{R}^N.$

The space $\mathcal{H}^\mathcal{L}(\mathbb{R}^N)$ is a Hilbert space; it is  equipped with the norm
\[
\|f\|^2_{\mathcal{H}^\mathcal{L}(\mathbb{R}^N)}=\int\limits_{\mathbb{R}^N}|a(\xi)\hat{f}(\xi)|^2d\xi.
\]
\begin{definition} The generalised solution of equation \eqref{1.1} in $\mathbb{R}^N$ is a function $ u\in C\left([0,\infty);L^2(\mathbb{R}^N)\right),$ such that $x^{-2\beta}\mathcal{D}_x^{2\alpha } u, \mathcal{L}u\in C\left((0,\infty);L^2(\mathbb{R}^N)\right).$\end{definition}
\begin{theorem}\label{th2}Let $\phi \in \mathcal{H}^\mathcal{L}(\mathbb{R}^N).$ Then the generalized solution of equation \eqref{1.1} satisfying conditions
\begin{equation}\label{1.2*}
u(0,y)=\phi(y),\,y\in\mathbb{R}^N,
\end{equation} and
\begin{equation}\label{1.3*}
\lim\limits_{x\rightarrow+\infty}u(x,y)\,\,\,\,\text{is bounded for almost every}\,\,\,\,y\in\mathbb{R}^N,
\end{equation}
exists, it is unique and can be represented as \begin{equation}\label{1.5*}u\left( {x,y} \right) = \int\limits_{\mathbb{R}^N}e^{-iy\xi} \hat{\phi}(\xi) E_{\alpha, 1+\frac{\beta}{\alpha}, \frac{\beta}{\alpha}} \left( {- \sqrt{a(\xi)} x^{\alpha+\beta}  }\right) d\xi,\,(x,y)\in \mathbb{R}_+\times\mathbb{R}^N,\end{equation} where $\hat{\phi}(\xi)=\frac{1}{(2\pi)^N}\int\limits_{\mathbb{R}^N}e^{-i\xi s}\phi(s)ds.$

In addition, the solution $u$ satisfies the following estimates:
\begin{align*}\|u\|_{C(\mathbb{R}_+;L^2(\mathbb{R}^N))}\leq\|\phi\|_{L^2(\mathbb{R}^N)},\end{align*}
\begin{align*}\sup\limits_{x\in(0,\infty)}\left\|x^{-2\beta}\mathcal{D}_x^{2\alpha}u\left( {x, \cdot} \right)\right\|_{L^2(\mathbb{R}^N)} \leq \|\phi\|_{\mathcal{H}^\mathcal{L}(\mathbb{R}^N)},\end{align*} and \begin{align*}\sup\limits_{x\in(0,\infty)}\|\mathcal{L}u\left({x, \cdot} \right)\|_{L^2(\mathbb{R}^N)}\leq\|\phi\|_{\mathcal{H}^\mathcal{L}(\mathbb{R}^N)}.\end{align*}
\end{theorem}

\subsection{Particular cases} We now specify Theorem \ref{th2} to several concrete cases.
\subsubsection{Laplace equation in the half-space}
Our first example will focus on the Laplace equation.

Let $\alpha=1,$ $\beta=0$ and $\mathcal{L}=-\Delta=\sum\limits_{j=1}^N\frac{\partial^2}{\partial y_j^2}.$ Then the equation \eqref{1.1} coincides with the classical Laplace equation on the half-space
\begin{equation}\label{1.1a*} u_{xx}(x,y) + \Delta_y u(x,y) = 0,\,\left({x,y}
\right) \in \mathbb{R}_+\times\mathbb{R}^N.\end{equation} It is known that the unique solution to problem \eqref{1.1a*}, \eqref{1.2*}, \eqref{1.3*} is represented by the Poisson integral (\cite{Stein})
\begin{equation*}u\left( {x,y} \right) = \frac{\Gamma((N+1)/2)}{\pi^{(N+1)/2}}\int\limits_{\mathbb{R}^N} \frac{x\phi(s)}{(|y-s|^2+x^2)^{(N+1)/2}}ds.\end{equation*}

\subsubsection{Multidimensional degenerate elliptic equations}
Let $\alpha=1,$ $\beta>-2$ and $\mathcal{L}=-\Delta_y.$

$\bullet$ If $\beta=1,$ then the equation \eqref{1.1} coincides with the multidimensional Tricomi equation \begin{equation}\label{1.1c*}u_{xx}(x,y)+x \Delta_y u(x,y)=0,\,x>0,\,y\in\mathbb{R}^N,\end{equation}
and the solution to problem \eqref{1.1c*}, \eqref{1.2*}, \eqref{1.3*} can be written as (\cite{Alg})
\begin{equation*}u\left( {x,y} \right) = \frac{3^{n+1/2}\Gamma(2/3)\Gamma(N/2+1/3)}{2^{1/3}\pi^{N/2+1}}\int\limits_{\mathbb{R}^N} \frac{x\phi(s)}{(9|y-s|^2+4x^3)^{N/2+1/3}}ds.\end{equation*}

$\bullet$ If $\beta=m>-2$ then the equation \eqref{1.1} coincides with the multidimensional Gellerstedt equation \begin{equation}\label{1.1cc*}u_{xx}(x,y)+x^m \Delta_y u(x,y)=0,\,x>0,\,y\in\mathbb{R}^N,\end{equation}
and the unique solution to problem \eqref{1.1cc*}, \eqref{1.2*}, \eqref{1.3*} can be written as (\cite{Alg})
\begin{equation*}u\left( {x,y} \right) = \frac{(m+2)^{n+\frac{1}{2}}\Gamma\left(\frac{2}{3}\right)\Gamma\left(\frac{N}{2}+\frac{1}{m+2}\right)}{2^{N}\pi^{\frac{N}{2}} \Gamma\left(\frac{1}{m+2}\right)}\int\limits_{\mathbb{R}^N} \frac{x\phi(s)}{\left(x^{m+2}+\left(\frac{m+2}{2}\right)^2|y-s|^2\right)^{\frac{N}{2}+\frac{1}{m+2}}}ds.\end{equation*}

\subsubsection{Fractional Laplace equation}
Let $\beta=0$ and $$\mathcal{L}v=(-\Delta)^sv=C_{N, s} P.V. \int_{\mathbb{R}^N}\frac{(v(y)-v(s))}{|x-y|^{N+2s}}dy,$$ where $s\in(0,1)$ and $C_{N,s}$ is a normalizing constant (whose value is not important here).
Then the equation \eqref{1.1} coincides with the equation \begin{equation}\label{1.1d*}\mathcal{D}_x^{2\alpha}u(x,y)+(-\Delta)^s_y u(x,y)=0,\,x>0,\,y\in\mathbb{R}^N.\end{equation}
From Theorem \ref{th2} we have the unique solution of the problem \eqref{1.1d*}, \eqref{1.2*}, \eqref{1.3*} in the form
\begin{equation*}u\left( {x,y} \right) = \int\limits_{\mathbb{R}^N}e^{-iy\xi} \hat{\phi}(\xi) E_{\alpha, 1} \left( {- |\xi|^s} x^{\alpha}  \right) d\xi,\,(x,y)\in \mathbb{R}_+\times\mathbb{R}^N.\end{equation*}
Rearranging the order of integration in the last representation, according to Fubini's Theorem, we have
\begin{equation*}u\left( {x,y} \right) = \int\limits_{\mathbb{R}^N}\phi(s) \int\limits_{\mathbb{R}^N} e^{-i\xi(y-s)} E_{\alpha, 1} \left( {- |\xi|^s} x^{\alpha}  \right) d\xi ds,\,(x,y)\in \mathbb{R}_+\times\mathbb{R}^N.\end{equation*}
Using the calculation of the Fourier transform of Mittag-Leffler functions from \cite{Zacher}, we have
\begin{equation*}u\left( {x,y} \right) = \pi^{-\frac{N}{2}}\int\limits_{\mathbb{R}^N}\frac{\phi(s)}{|y-s|^{N}} H_{3\,2}^{1\,2}\left(\frac{2^s x^\alpha}{|y-s|^{s}}\Big{|} \begin{array}{l}(1-N/2-s/2), \,\,(0,1),\,\, (0,s/2)\\ (0,1),\,\, (0,\alpha) \end{array}\right) ds.\end{equation*}
Here $H_{pq}^{mn}(\cdot)$ is the Fox H-function defined via a Mellin-Barnes type integral as
$$H_{p\,q}^{m\,n}\left(z\Big{|} \begin{array}{l}(a^1_i, a^2_i)_{1,p}\\ (b^1_i, b^2_i)_{1,q}\end{array}\right)=\frac{1}{2\pi i}\int\limits_{\mathcal{I}}\mathcal{H}_{p, q}^{m, n}(\tau) z^{-\tau} d\tau,$$ where $(a^1_i, a^2_i)_{1,p}=((a^1_1, a^2_1), (a^1_2, a^2_2), ..., (a^1_p, a^2_p))$ and $$\mathcal{H}_{p, q}^{m, n}(\tau)=\frac{\prod\limits_{j=1}^m\Gamma(b^1_j+b^2_j\tau) \prod\limits_{i=1}^n\Gamma(1-a^1_i-a^2_j\tau)}{\prod\limits_{i=n+1}^p\Gamma(a^1_i+a^2_i\tau)\prod\limits_{j=m+1}^q\Gamma(1-b^1_j-b^2_j\tau)}.$$

\subsection{Proof of Theorem \ref{th2}}
\subsubsection{Existence of solution.} Applying the Fourier transform $\mathcal{F}$ to problem \eqref{1.1}, \eqref{1.2*}--\eqref{1.3*} with respect to space variable $y$ yields
\begin{equation}\label{3.3*}\mathcal{D}_x^{2\alpha } \hat{u}
\left( x, \xi\right) - a(\xi) x^{2\beta}\hat{u}\left( x, \xi \right) = 0,\, x >0,\,\xi\in \mathbb{R}^N,\end{equation} \begin{equation}\label{3.4*}\hat{u} \left(0,\xi\right) = \hat{\phi}(\xi),\, \hat{u} \left( \infty, \xi\right)
\,\,\text{is bounded for}\,\,\xi\in \mathbb{R}^N,\end{equation} thank to $\mathcal{F}\left\{\mathcal{L} u(x,y)\right\}=a(\xi)\hat{u}(x,\xi).$ Then the solution of problem \eqref{3.3*}-\eqref{3.4*} can be represented as
\begin{equation}\label{3.6*}\hat{u} \left(
x, \xi\right) = \hat{\phi}(\xi) E_{\alpha, 1+\frac{\beta}{\alpha}, \frac{\beta}{\alpha}} \left( {- \sqrt{a(\xi)} x^{\alpha+\beta}  }\right).\end{equation}
By applying the inverse Fourier transform $\mathcal{F}^{-1}$ we have \eqref{1.5*}, i.e.
\begin{equation*}u\left( {x,y} \right) = \int\limits_{\mathbb{R}^N}e^{iy\xi} \hat{\phi}(\xi) E_{\alpha, 1+\frac{\beta}{\alpha}, \frac{\beta}{\alpha}} \left( {- \sqrt{a(\xi)} x^{\alpha+\beta}  }\right) d\xi,\,(x,y)\in \mathbb{R}_+\times\mathbb{R}^N.\end{equation*}

\subsubsection{Convergence of solution.} Now we prove the convergence of the obtained solution.
Applying estimate \eqref{2.6} and Plancherel theorem we have
\begin{align*}\sup\limits_{x\in[0,\infty)}\int\limits_{\mathbb{R}^N}|u\left( {x,y} \right)|^2dy&=\sup\limits_{x\in[0,\infty)}\int\limits_{\mathbb{R}^N}|\hat{u}\left( {x,\xi} \right)|^2d\xi\\& \leq  \sup\limits_{x\in[0,\infty)}\int\limits_{\mathbb{R}^N}\left|\hat{\phi}(\xi)\right|^2 \left|E_{\alpha, 1+\frac{\beta}{\alpha}, \frac{\beta}{\alpha}} \left( {- \sqrt{a(\xi)} x^{\alpha+\beta}  }\right)\right|^2d\xi\\& \leq \int\limits_{\mathbb{R}^N}\left|\hat{\phi}(\xi)\right|^2 d\xi=\|\hat{\phi}\|^2_{L^2(\mathbb{R}^N)}=\|{\phi}\|^2_{L^2(\mathbb{R}^N)}<\infty.\end{align*}
Let us calculate $\mathcal{D}_x^{2\alpha } u:$
\begin{align*}\mathcal{D}_x^{2\alpha}u\left( {x,y} \right)& = \int\limits_{\mathbb{R}^N}e^{iy\xi} \hat{\phi}(\xi) \mathcal{D}_x^{2\alpha}E_{\alpha, 1+\frac{\beta}{\alpha}, \frac{\beta}{\alpha}} \left( {- \sqrt{a(\xi)} x^{\alpha+\beta}  }\right) d\xi\\&= x^{2\beta}\int\limits_{\mathbb{R}^N}e^{iy\xi} \hat{\phi}(\xi)a(\xi)E_{\alpha, 1+\frac{\beta}{\alpha}, \frac{\beta}{\alpha}} \left( {- \sqrt{a(\xi)} x^{\alpha+\beta}  }\right) d\xi,\,\,(x,y)\in \mathbb{R}_+\times\mathbb{R}^N.\end{align*}
Hence
\begin{align*}\sup\limits_{x\in(0,\infty)}\left\|x^{-2\beta}\mathcal{D}_x^{2\alpha}u\left( {x, \cdot} \right)\right\|^2_{L^2(\mathbb{R}^N)} &\leq \sup\limits_{x\in(0,\infty)}\int\limits_{\mathbb{R}^N}a^2(\xi)|\hat{\phi}(\xi)|^2\left|E_{\alpha, 1+\frac{\beta}{\alpha}, \frac{\beta}{\alpha}} \left( {- \sqrt{a(\xi)} x^{\alpha+\beta}  }\right)\right|^2 d\xi\\& \leq \int\limits_{\mathbb{R}^N}|a(\xi)\hat{\phi}(\xi)|^2d\xi=\|\phi\|_{\mathcal{H}^\mathcal{L}(\mathbb{R}^N)}^2<\infty.\end{align*} Similarly, for $\mathcal{L}u$ we have
\begin{align*}\sup\limits_{x\in(0,\infty)}\left\|\mathcal{L}u\left( {x, \cdot} \right)\right\|^2_{L^2(\mathbb{R}^N)} &\leq \sup\limits_{x\in(0,\infty)}\int\limits_{\mathbb{R}^N}a^2(\xi)|\hat{\phi}(\xi)|^2\left|E_{\alpha, 1+\frac{\beta}{\alpha}, \frac{\beta}{\alpha}} \left( {- \sqrt{a(\xi)} x^{\alpha+\beta}  }\right)\right|^2 d\xi\\& \leq \|\phi\|_{\mathcal{H}^\mathcal{L}(\mathbb{R}^N)}^2<\infty.\end{align*}

\subsubsection{Uniqueness of solution.} Suppose that there are two solutions $u_1(x, y)$ and $u_2(x,y)$ of problem \eqref{1.1}, \eqref{1.2*}--\eqref{1.3*}. Let $u(x,y)=u_1(x,y)-u_2(x,y).$ Then $u(x,y)$ satisfies the equation \eqref{1.1} and homogeneous conditions \eqref{1.2*}--\eqref{1.3*}.

Let us consider the function \begin{equation}\label{3.5*}\hat{u}(x,\xi)=\int\limits_{\mathbb{R}^N}e^{-iy\xi} u(x,y)dy,\,\,x\geq 0,\,\xi\in\mathbb{R}^N.\end{equation}
As $u$ is bounded continuous in $x$ function, applying $\mathcal{D}^{2\alpha}_x$ to the function \eqref{3.5*} by \eqref{1.1} we have
\begin{align*}\mathcal{D}^{2\alpha}_x\hat{u}(x,\xi)&=\int\limits_{\mathbb{R}^N}e^{-iy\xi}\mathcal{D}^{2\alpha}_x u(x,y)dy\\&=x^{2\beta}\int\limits_{\mathbb{R}^N}e^{-iy\xi}\mathcal{L}u(x,y)dy \\& =x^{2\beta}\mathcal{F}\left[\mathcal{F}^{-1}(a(\xi)\hat{u}(x,y))\right]\\& = x^{2\beta}a(\xi) \hat{u}(x,\xi),\,x\geq 0,\, \xi\in\mathbb{R}^N.\end{align*} Also from \eqref{1.2*} and \eqref{1.3*} we have
$\hat{u}(0,\xi)=0,\,\,\hat{u}(\infty,\xi)\,\,\text{is bounded}.$ Then from \eqref{3.6*} we conclude that $\hat{u}(x,\xi)=0,\,x\geq 0,\,\xi\in\mathbb{R}^N.$ Applying the inverse Fourier transform we have $u(x,y)\equiv 0,\,(x,y)\in [0,\infty)\times\mathbb{R}^N.$
The proof is complete.

\end{document}